\documentclass[12pt,a4paper,leqno]{article}
\usepackage{amsmath}
\usepackage{amssymb}
\usepackage{amsfonts}
\usepackage{amsthm}
\usepackage{xypic}
\usepackage[dvips]{graphicx,epsfig,color}
\usepackage[width=15cm]{geometry}
\usepackage[english]{babel}
\usepackage{mathrsfs}
\usepackage{lscape}
\usepackage{longtable}
\usepackage{xtab}
\usepackage[applemac]{inputenc}
\numberwithin{equation}{section}
\renewcommand{\phi}{\varphi}

\DeclareMathOperator{\Image}{Im}

\DeclareMathOperator{\rank}{rank}

\def \Z {{\mathbb Z}}
\def \C {\mathbb C}
\def \R {\mathbb R}
\def \cH {\mathcal H}
\def \cM {\mathcal M}
\def \tcM {\widetilde{\mathcal M}}

\def \cU {\mathcal U}

\def \ev {{\mathrm{ev}}}

\newcommand{\p}{\partial}

\newtheorem{theorem}{Theorem}[section]
\newtheorem{lemma}[theorem]{Lemma}
\newtheorem{proposition}[theorem]{Proposition}
\newtheorem{Cor}[theorem]{Corollary}

\theoremstyle{definition}
\newtheorem{definition}[theorem]{Definition}
\newtheorem{example}[theorem]{Example}
\newtheorem{examples}[theorem]{Examples}

  \newtheorem*{def1.1}{Definition \ref{C0-comm}}

\theoremstyle{remark}
\newtheorem{remark}[theorem]{Remark}

\title{\textbf {On the topology of fillings of contact manifolds and applications}\thanks{Supported by ANR project "Floer Power" ANR-08-BLAN-0291-03/04. 
We are grateful to Beijing International Center for Mathematical Research of Beijing University 
for hospitality during the completion of this paper.}
}
\author{  \\{ \sc Alexandru Oancea${\ }^{1}$\qquad Claude Viterbo${\ }^{2}$}\\
\\
${\ }^{1}$ IRMA \\ UMR 7501, Universit\'e de Strasbourg \\ 7 rue Ren\'e
Descartes, 67084 Strasbourg \\ {\tt   oancea@math.u-strasbg.fr}
\\
\\${\ }^{2}$  Centre de Math\'ematiques Laurent Schwartz \\ UMR 7640
du CNRS \\ \'Ecole Polytechnique - 91128 Palaiseau, France
\\ {\tt viterbo@math.polytechnique.fr   } }
\date{November 28, 2009}

\begin{document}

\maketitle

\begin{abstract}The aim of this paper is to address the following
question: given a contact manifold $(\Sigma, \xi)$, what can be said about the
symplectically aspherical manifolds $(W, \omega)$ bounded by $(\Sigma, \xi)$ ?
We first extend a theorem of Eliashberg, Floer and McDuff to prove that, under suitable assumptions, the map from $H_{*}(\Sigma)$ to $H_{*}(W)$ induced by inclusion is surjective.
We apply this method in the case of contact manifolds admitting a contact embedding in ${\mathbb R}^{2n}$ or in a subcritical Stein manifold. We prove in many cases that the homology of the fillings is uniquely determined.  
Finally, we use more recent methods of symplectic topology to prove that,   
if a contact hypersurface has a subcritical Stein filling, then all its $SAWC$ fillings have the same homology. 

A number of applications are given, from obstructions to the existence of Lagrangian or contact embeddings, to the  exotic nature of some contact structures. We refer to the table in Section~\ref{table} for a summary of our results. 
\end{abstract}

\section{Introduction}
In this paper all symplectic  manifolds will be assumed to be connected, of dimension $2n$, and symplectically aspherical, meaning that the symplectic form vanishes on the second homotopy group. All contact manifolds are connected and have dimension $2n-1$. We denote by $\sigma _{0}$ the standard symplectic form on 
${\mathbb R}^{2n}$ or ${\mathbb C}P^n$, and by $\alpha_{0}$ the standard contact form on the sphere $S^{2n-1}$.

In a celebrated paper (\cite{McDuff}), Eliashberg, Floer and McDuff proved that, if $(W, \omega)$ is a symplectically aspherical manifold with contact 
boundary $(S^{2n-1}, \alpha _{0})$, then $W$ is diffeomorphic  to the unit ball $B^{2n}$. In the case of dimension $4$, Gromov had 
earlier proved (\cite{Gromov}) that $W$ is actually symplectomorphic to $(B^{4}, \sigma_{0})$, but this relies heavily on positivity of 
intersection for holomorphic curves that is special to dimension $4$.

One can ask more generally, given a fillable contact manifold  $(\Sigma, \xi)$ and a symplectically aspherical filling $(W, \omega)$, what can be 
said about the topology or the homology of $W$. Is it  uniquely determined by the contact structure $(\Sigma, \xi)$~? Is it 
determined by the topology of $\Sigma$~? Do we have lower bounds~?  Upper bounds~?
 It turns out that all these possibilities actually occur.

For example, if $(\Sigma, \xi)$ has a contact embedding into $( {\mathbb R}^{2n}, \sigma_{0})$ - many such examples can be found 
in \cite{Laudenbach} - it readily follows from the Eliashberg-Floer-McDuff theorem and some elementary algebraic topology that all subcritical Stein fillings have the same homology. If the homology of $\Sigma$ vanishes in degree $n$, we can prove that all Stein fillings have the same homology. This gives easy examples of contact manifolds with no 
contact embedding in $({\mathbb R}^{2n}, \sigma_{0})$. As far as the authors know, there are only  
few previously known examples of fillable manifolds not embeddable in ${\mathbb R}^{2n}$, with the 
exception of recent results in \cite{Cieliebak-Frauenfelder-Oancea} and \cite{Albers-McLean}, which however assume
the exactness of the embeddings, an assumption we usually can dispense with. 
More general results of the same homological flavour follow from the same methods, and a generalization of the Eliashberg-Floer-McDuff theorem 
to the case of subcritical Stein manifolds. These are manifolds $W$ admitting an 
exhausting plurisubharmonic function with no critical points of index $n=\frac{1}{2} \dim (W)$ (see Definition~\ref{def:Stein}).

Our last result  uses more sophisticated tools. One of them will be symplectic homology of $W$, and its positive part, 
defined in \cite{Viterbo1}. It turns out that this positive part, under mild assumptions on the Conley-Zehnder index of closed characteristics, 
is independent of the filling. This is proved in \cite{Cieliebak-Frauenfelder-Oancea} as a consequence of arguments in \cite{Bourgeois-Oancea-2}.  A symplectically aspherical manifold $(W,\omega)$ with contact type boundary is called a \emph{$SAWC$-manifold} if its symplectic homology vanishes (this is equivalent to the \emph{Strong Algebraic Weinstein Conjecture (SAWC)} formulated in \cite{Viterbo1}, cf. Section~\ref{shof}).  We show that, if $(\Sigma , \xi )$ bounds a subcritical Stein manifold $(W, \omega)$, any other $SAWC$ filling will have the same homology as $W$. 

Of course many questions remain open.  As far as we can see, nothing can be said about the symplectic topology of fillings 
outside the subcriticality/non-subcriticality alternative. Are there examples of compact manifolds $L$ such that $ST^*L$ has fillings with homology different from $H_{*}(L)$ ? Is there an embedding of the Brieskorn sphere of a singularity of Milnor number $\mu$ in the Milnor fibre of a singularity of Milnor number $\mu'<\mu$ ? 

\medskip 

\noindent {\it Acknowledgements.} We are grateful to Vincent Blanl\oe il for his help with the formulation and proof of Proposition~\ref{prop:Brieskorn}.

\section{The Eliashberg-Floer-McDuff theorem\\revisited}

\noindent {\bf Conventions.} We denote by $(W,\omega)$ a symplectic manifold of dimension $2n$ which is symplectically aspherical ($[\omega]\pi_2(W)=0$). We denote by $(\Sigma,\xi)$ a contact manifold of dimension $2n-1$.  
We assume that $\xi$ is co-orientable, and fix 
a co-orientation. The contact structure $\xi$ is then defined by a contact form $\alpha$, and $\Sigma$ is oriented by   
$\alpha \wedge (d\alpha)^{n-1}\neq 0$.
All homology and cohomology groups are taken with coefficients in a field. 

\begin{definition}
A {\it contact embedding} of $(\Sigma , \xi)$ in $(W, \omega)$ is a codimension $1$  embedding  such that there exists a positive contact form $\alpha$ extending to a neighbourhood of $\Sigma$ as a primitive of $\omega$. The contact embedding is called {\it exact} if $\alpha$ extends to the whole of $W$ as a primitive of $\omega$. 
\end{definition}

\begin{definition}
A (co-oriented) hypersurface $\Sigma\subset (W,\omega)$ is said to be of {\it contact type} in $W$ if there exists a primitive $\alpha$ of $\omega$, defined in a neighbourhood of $\Sigma$, and restricting on $\Sigma$ to a contact form (whose $\omega$-dual vector field defines the positive co-orientation of $\Sigma$). 
The hypersurface is said to be of {\it restricted contact type} in $W$ if there exists such a primitive $\alpha$ which is globally defined on $W$.
\end{definition}
\begin{definition}
A {\it symplectic filling} of $(\Sigma , \xi)$ is a symplectic manifold $(W,\omega)$ 
without closed components, such that $\partial W=\Sigma$ and there exists a positive contact form $\alpha$ extending to a 
neighbourhood of $\Sigma$ as a primitive of $\omega$. We shall say that the symplectic filling is {\it exact} 
if $\alpha$ extends to the whole of $W$ as a primitive of $\omega$.
\end{definition}

 \begin{definition} \label{def:Stein}
 A symplectic filling $(W,\omega)$ of $(\Sigma, \xi)$ is a {\it Stein filling}  if
 $W$ has a complex structure $J$, and a non-positive plurisubharmonic function $\psi$, such that
 $\Sigma= \psi^{-1}(0)$ and $-J^*d\psi$ is a contact form defining $\xi$.
 Note that  $\psi$ can always be chosen to be a Morse function. Then its critical points have index at most $n$, so that $W$ has the homotopy type of a CW complex of dimension $\le n$. If we can find the function $\psi$ with no critical points of index $n$, then $W$ is said to be {\it subcritical Stein}.
 \end{definition}

\begin{remark} A contact embedding of $(\Sigma, \xi)$ in $(W, \omega)$ which is {\it separating} -- i.e. $W\setminus \Sigma$ consists of two connected components --
yields a filling of $(\Sigma, \xi)$ by the connected component of $W\setminus \Sigma$ for which the boundary orientation of $\Sigma$ coincides 
with the orientation induced by a positive contact form $\alpha$. This filling we shall call {\it the interior of} $\Sigma$ and we shall denote by $Z$. If $W$ is non-compact, $Z$ is the bounded component of $W\setminus \Sigma$. Note that $\Sigma$ is always separating if $H_{2n-1}(W;{\mathbb Z})$ is torsion, and in particular if $W$ is Stein and $n\ge 2$.
 \end{remark}

Our goal in this section is  to prove the following theorem.

\begin{theorem}  \label{prop:steinsubcritfillings}
Assume $(\Sigma, \xi )$ admits a contact embedding in a subcritical Stein manifold $(M,\omega_0)$, with interior component $Z$. 
Let $(W, \omega)$ be a symplectically aspherical filling of $\Sigma$ and assume that one of the following two conditions is satisfied:
\begin{enumerate} 
\item \label{a} $H_{2}(W,\Sigma)=0$;
\item \label{b} $\Sigma$ is simply connected.
\end{enumerate} 
Then the map 
$$
H_{j}(\Sigma) \longrightarrow H_{j}(W)
$$ 
induced by inclusion 
is onto in every degree $j\ge 0$. 
\end{theorem}

\begin{remark}
 Condition \ref{a} holds if $W$ is Stein. The embedding $\Sigma\hookrightarrow M$ is always separating since $H_{2n-1}(M;\Z)=0$.
 \end{remark} 
 \begin{remark} 
When $\Sigma$ is a sphere, we get that $W$ has vanishing homology. This is the original Eliashberg-Floer-McDuff theorem 
(see \cite{McDuff}), since an application of the $h$-cobordism theorem, plus the fact - due to Eliashberg- that $\pi_{1}(W)$ vanishes,
implies that $W$ is diffeomorphic to the ball $B^{2n}$. Indeed, since $H_{j}(S^{2n-1})=0$ for $1\leq j \leq 2n-2$, the same holds for $H_{j}(W)$. In particular $H_1(W)=0$, which implies that $H^1(W)$ and $H_{2n-1}(W)$ vanish. When $n=2$ Gromov (see \cite{Gromov}) proved that 
$W$ is symplectomorphic to the ball $B^4$, but this relies heavily on purely $4$-dimensional arguments 
(positivity of intersection of holomorphic curves).  
 \end{remark}

 Our proof of Theorem~\ref{prop:steinsubcritfillings} closely follows the original proof in \cite{McDuff}, 
except for the final homological argument. 

We shall start by working in the following special situation, and we will then prove that this is enough to deal with the  general case.

Let  $(P,\omega_{P})$ be a symplectic manifold and $H$ be a codimension two symplectic submanifold of $P$. 
  
\begin{figure}[htp] 
\centering  
\includegraphics[width=14cm]{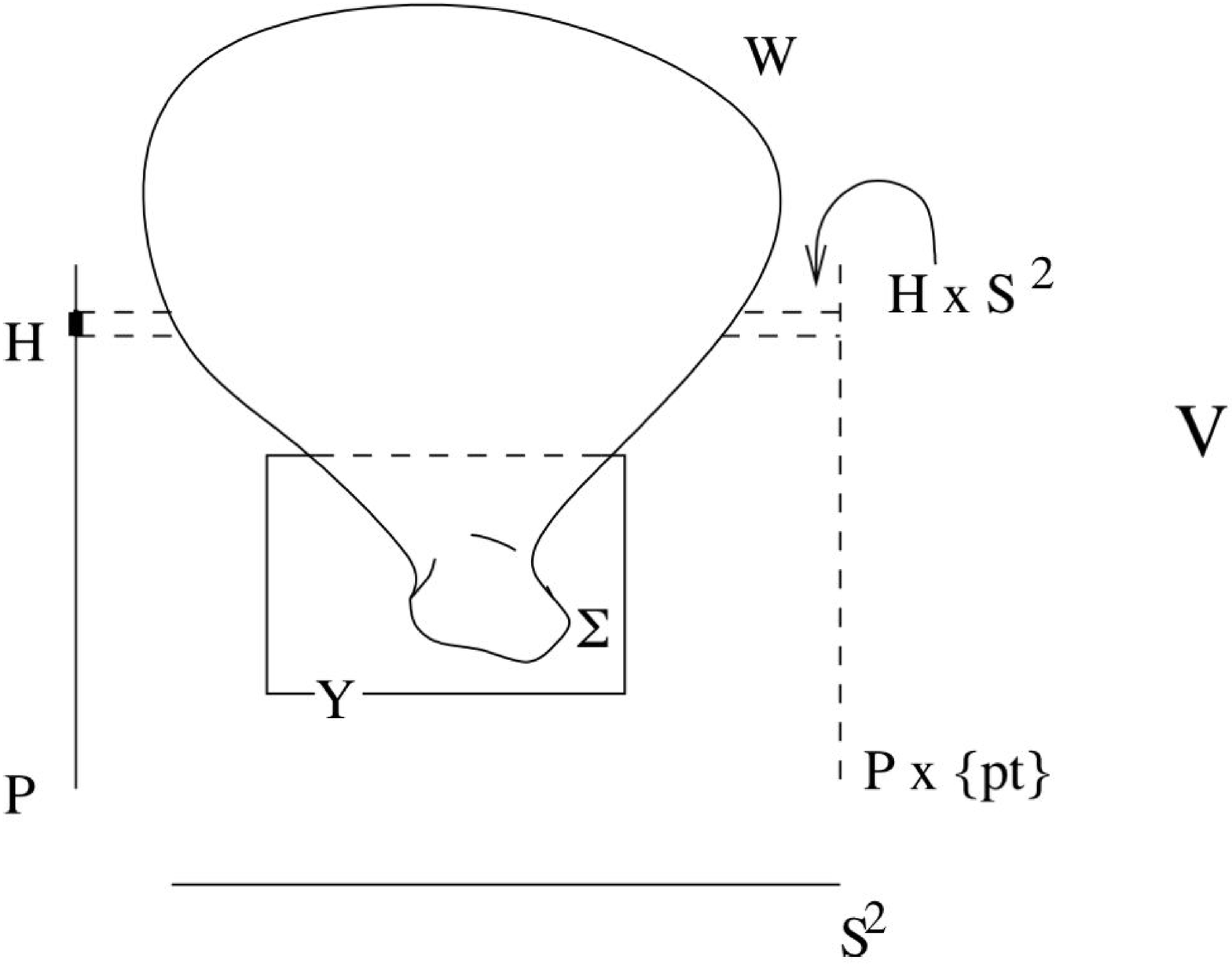} 
\caption{{Geometric setup}}\label{fig:greffe} 
\end{figure}

We consider the symplectic manifold $(P\times S^2, \omega_{P}\oplus \sigma)$, where $\sigma$ is the standard symplectic form on $S^2$ normalized by $[\sigma][S^2]=1$. Viewing $S^2$ as $\C\cup\{\infty\}$ we denote $D_-^2:=\{z\, : \, |z|<\frac 1 2\}$ and $D_+^2:=\{z\, : \, |z|>\frac 1 2\}$. Let $(\Sigma, \alpha)$ be a separating contact manifold  contained in $(P\setminus H)\times D_{-}^2$, with interior $Z$. We set $Y = (P\times S^2 \setminus Z)$ and
 $$
 V =Y \sqcup_{\Sigma} W= (P\times S^2 \setminus Z) \sqcup_{\Sigma} W,
 $$ 
where $(W,\omega)$ is a filling of $(\Sigma , \alpha)$ (Figure~\ref{fig:greffe}). Then $V$ has a symplectic form $\omega_{V}$ obtained by gluing $\omega_{P}\oplus \sigma$ on $Y$ and $\omega$ on $W$. 

Let  $p_{0}$ be a point in $H$ and denote $A:=[\{p_0\}\times S^2]\in H_2(V;\Z)$ (note that, for $p_{0}\in H$,  we have that $\{p_{0}\}\times S^2 \subset Y\subset V$). Given an $\omega_V$-compatible almost complex structure $J$ on $V$, we denote by $\tcM_J$ the space of $J$-holomorphic maps $u: {\mathbb C}P^1  \longrightarrow V$ representing the class $A$. 

\begin{lemma} \label{lem:simple}
Assume the symplectic form $\omega_P$ is integral, i.e. $[\omega_P]\in\mathrm{Im}(H^2(P;\Z)\to H^2(P;\R))$. Assume also that one of the following conditions holds: 
\begin{enumerate}
\item \label{item:lem-1} $H_2(W,\Sigma)=0$;
\item \label{item:lem-2} $\Sigma$ is simply connected. 
\end{enumerate}
Then, for any $\omega_V$-compatible almost complex structure $J$ on $V$, the class $A$ is $J$-simple, meaning that it cannot be decomposed as $A=B+C$, with $B,C\in H_2(V;\Z)$ represented by  non-constant $J$-holomorphic spheres.  
\end{lemma}

\begin{proof} Arguing by contradiction, let us assume that such a decomposition exists for some $\omega_V$-compatible $J$. We obtain in particular $0<\omega_V(B),\omega_V(C)<\omega_V(A)$. Let 
$\Gamma_B$ and $\Gamma_C$ be the  $J$-holomorphic spheres representing $B$ and $C$. Since $W$ is symplectically aspherical, it is not possible that any of the cycles $\Gamma_B$, $\Gamma_C$ be entirely contained in $W$. Since $\omega_P$ is integral and $\int_{S^2}\sigma=1$, it follows that the class $A$ has minimal area in $P\times S^2$, so that it is neither possible that any of the cycles $\Gamma_B$, $\Gamma_C$ be entirely contained in $Y\subset P\times S^2$. 
Thus $\Gamma_B$ and $\Gamma_C$ intersect both $Y$ and $W$. By (smoothly) perturbing the representing $J$-holomorphic spheres, we can achieve transverse intersection with $\Sigma$, along some collection of circles. 

Let us now assume~\ref{item:lem-1}. We consider the two pieces in $\Gamma_C$ separated by $\Sigma$. We denote by $C_{1}$ the piece in $W$ and 
by $C_{2}$ its complement. 
Then $\int_{\Gamma_C}\omega_V= \int_{C_{1}}\omega_V + \int_{C_{2}}\omega_V$. But $C_{1} \in H_{2}(W,\Sigma)$ so by our assumption there is a cycle 
$\Gamma$ in $\Sigma$ such that $\partial C_{1}=\partial \Gamma$. Because $C_{1}\cup \Gamma$ is a cycle in $H_{2}(W)$, and the map 
$H_{2}(\Sigma) \longrightarrow H_{2}(W)$ is onto (using again $H_{2}(W,\Sigma)=0$), we obtain that  $C_{1}\cup \Gamma$ is homologous 
in $W$ to a cycle $C_{3}$ contained in $\Sigma$. Thus $ \int_{C_{1}}\omega_V
-\int_{\Gamma}\omega_V = \int_{C_{3}}\omega_V$, which vanishes because $\omega_V$ is exact near $\Sigma$.  Finally $C_{2}\cup \Gamma$ is a 
cycle in $Y$ with the same area as $\Gamma_C=C_{1}\cup C_{2}$. Since $Y\subset P\times S^2$, this contradicts the fact that the class $A$ has minimal area in $P\times S^2$.

We now assume~\ref{item:lem-2}. As above, let $C_{1}, C_{2}$ be the parts of $\Gamma_C$ separated by $\Sigma$, with $C_1$ being the piece contained in $W$. By assumption, we can cap all the common boundary circles of $C_1$ and $C_2$ by discs. Let us denote this union of discs by $\Gamma$. Then
$C_1\cup \Gamma$ is a collection of spheres in $W$ and, by symplectic asphericity, it has zero area. Thus $C_2\cup\Gamma$ is a cycle in $Y$ with the same area as $\Gamma_C=C_{1}\cup C_{2}$. This contradicts again the minimality of the area of the class $A$ in $P\times S^2$.  
\end{proof}

As a consequence of Lemma~\ref{lem:simple} we have the following facts: 
\begin{itemize}
\item for any $\omega_V$-compatible almost complex structure $J$ on $V$, the elements of $\tcM_J$ are simple curves (i.e. they are not multiply covered);
\item for a generic choice of the $\omega_V$-compatible almost complex structure $J$, the linearized Cauchy-Riemann operator is surjective for every element of $\tcM_J$ and $\tcM_J/PSL(2,\C)$ is a smooth manifold of dimension 
$$
\dim\, \tcM_J/PSL(2,\C) = 2n+2\langle c_1(V),A\rangle - 6 = 2n-2.
$$
Such an almost complex structure $J$ is called \emph{regular}.
\item the manifold $\tcM_J/PSL(2,\C)$ is compact.
\end{itemize}

Let $J$ be a regular almost complex structure. It is convenient to consider the following model for the manifold $\tcM_J/PSL(2,\C)$. Given (disjoint) submanifolds $T_{-1},T_1,T_\infty\subset V$ that are $C^2$-close to $P\times\{-1\}$, $P\times \{1\}$, respectively $P\times\{\infty\}$, we denote by $\cM:=\cM_J$ the set of elements $u\in\tcM_J$ such that 
$u(z)\in T_z$, $z=-1,1,\infty$. For a generic choice of the perturbations $T_z$, $z=-1,1,\infty$ the relevant evaluation maps are transverse, so that $\cM$ is a compact submanifold of $\tcM_J$ of dimension $2n-2$. Since $\cM$ intersects every orbit of $PSL(2,\C)$ exactly once, it follows that the natural map $\cM\to \tcM/PSL(2,\C)$ is a diffeomorphism. 

Given an $\omega_P$-compatible almost complex structure $J_P$ on $P$, we denote by $\widetilde J_P$ the almost complex structure $J_P\oplus i$ on $P\times S^2$. 
For the proof of the next result, we closely follow \cite{McDuff}, pp. 660--661.

Let $J_P$ be an $\omega_P$-compatible almost complex structure on $P$. We say that $H$ is {\it hyperplane section-like for $J_P$} if the following hold: 
\begin{itemize}
\item $H$ is a $J_P$-complex submanifold,
\item there exists a codimension two $J_P$-complex (singular) submanifold $B\subset H$ (the {\it base locus}),
  a relatively compact neighbourhood $\cU$ of $H$ and a relatively compact open neighbourhood $\cU_B$ of $B$,
and a family $\cH_z$ of $J_P$-complex hypersurfaces parametrized by an open neighbourhood of $0$ in $\C$ and contained in $\cU$,
such that $\cH_z\cap (\cU\setminus \overline\cU_B)$ foliate $\cU\setminus \overline\cU_B$.
\end{itemize}

Note that the hyperplane section of a projective manifold is hyperplane-section-like for the underlying complex structure. 

\begin{lemma} \label{lem:degree} Let $J_P$ be an $\omega_P$-compatible almost complex structure on $P$ such that 
$H$ is hyperplane section-like for $J_P$. For every regular $J$ which is close to $\widetilde J_P$ on a neighbourhood of $H\times S^2\subset V$, the evaluation map 
$$
\ev:\cM\times S^2\to V
$$
has degree $\pm 1$.
\end{lemma}

\begin{proof} Let $B$ be the base locus for $H$ and let $\cU$, $\cU_B$ be the neighbourhoods of $H$ and $B$ such that $\cU\setminus \overline \cU_B$ is foliated by $J_P$-complex hypersurfaces $\cH\cap (\cU\setminus \overline\cU_B)$. We consider $\cU\times S^2$ as a neighbourhood of $H\times S^2$ in $V$, so that $(\cU\setminus\overline \cU_B)\times S^2$ is also foliated by the $\widetilde J_P$-complex hypersurfaces $(\cH\times S^2)\cap (\cU\setminus\overline \cU_B)\times S^2$. We prove the lemma in three steps. 

\smallskip \noindent {\it Step~1: Let $J$ coincide with $\widetilde J_P$ on $\cU\times S^2$. For every $p\in(\cU\setminus\overline\cU_B)\times S^2$, there exists a unique element of $\tcM_J/PSL(2,\C)$ through $p$.}
\smallskip 

Let $\cH_p\times S^2$ be the $\widetilde J_P$-complex hypersurface through $p$. Given a curve $C$ through $p$ represented by some $[u]\in\tcM_J/PSL(2,\C)$, the homological intersection between $C$ and $\cH_p\times S^2$ is zero. Since $C\cap \cH_p\times S^2\neq \emptyset$, it follows from the positivity of intersections for holomorphic curves that $C$ is entirely contained in $\cH_p\times S^2$. But, in $\cH_p\times S^2$, there is clearly a unique $\widetilde J_P$-holomorphic curve in the class $\{pt\}\times S^2$ through each point. (Remark: Positivity of intersections is proved with all details in a $4$-dimensional context in~\cite{Lazzarini-McDuff-Salamon}. The higher dimensional case of a curve intersecting a complex hypersurface is treated using exactly the same methods.)

\smallskip \noindent {\it Step~2: Let $J_0$ be an almost complex structure which coincides with $\widetilde J_P$ on $\cU\times S^2$. For every $J$ which is close enough to $J_0$, and for every point $p\in(\cU\setminus\overline\cU_B)\times S^2$, there is a unique element of $\tcM_J/PSL(2,\C)$ through $p$.}
\smallskip

We follow \cite{McDuff}, Lemma~3.5. Arguing by contradiction, we find a point $p\in (\cU\setminus\overline \cU_B)\times S^2$ and a sequence $J_\nu$, $\nu\ge 1$ converging to $J_0$ such that, for every $\nu$, there exist two distinct unparametrized $J_\nu$-holomorphic spheres $C_\nu$ and $C'_\nu$ through $p$. Since $A$ is a simple class (Lemma~\ref{lem:simple}), both $C_\nu$ and $C'_\nu$ converge as unparametrized spheres to the unique $J_0$-holomorphic sphere $C$ through $p$. In particular, for $\nu$ large enough they are both contained in $\cU'\times S^2$, for some relatively compact open subset $\cU'\subset \cU$.

We now view $\cU\times S^2$ as a subset of $P\times S^2$ and extend $J_\nu|_{\cU'\times S^2}$ to an almost complex structure $J'_\nu$ on $P\times S^2$ which is compatible with $\omega_P\oplus \sigma$ and satisfies $J'_\nu\to \widetilde J_P$, $\nu\to \infty$. For $\nu$ large enough, the $J_\nu$-holomorphic curves $C_\nu$, $C'_\nu$ passing through $p$ are now viewed in $P\times S^2$, where they are $J'_\nu$-holomorphic. The almost complex structure $\widetilde J_P$ is obviously regular for curves in the class $[\{pt\}\times S^2]$ in $P\times S^2$, and Step~1 shows that the evaluation map
$$
\ev:\tcM_{\widetilde J_P}\times_{PSL(2,\C)} S^2\to P\times S^2
$$
is a diffeomorphism (we make a slight abuse of notation and write here $\tcM_J$ for the space of $J$-holomorphic curves in $P\times S^2$ representing the class $[\{pt\}\times S^2]$). The evaluation map remains a diffeomorphism for small perturbations of $\widetilde J_P$, and we reach a contradiction with the fact that $p$ has at least two preimages via the evaluation maps $\ev:\tcM_{J'_\nu}\times_{PSL(2,\C)} S^2\to P\times S^2$ for $\nu$ large enough. 

\smallskip \noindent {\it Step~3: We prove the Lemma.}
\smallskip

The degree of the evaluation map can be computed by counting the number of preimages of a generic point in $V$. We can therefore choose our point generically in $(\cU\setminus\overline\cU_B)\times S^2\subset V$, and the number of preimages is then equal to one by Step~2.
\end{proof}

\begin{proposition}  \label{thm:BigMac}
Let $(P,\omega_P)$ be a symplectic manifold such that $\omega_P$ is an integral class, and $H\subset P$ be a codimension two symplectic submanifold which is hyperplane section-like for some $\omega_P$-compatible almost complex structure $J_P$. Let $(\Sigma,\xi)$ be a contact separating hypersurface of $(P\times S^2,\omega_P\oplus\sigma)$ which is contained in $(P\setminus H)\times S^2$. Let $W$ be any symplectically aspherical filling of $(\Sigma,\xi)$. Assume one of the following two conditions holds: 
\begin{enumerate}
\item $H_2(W,\Sigma)=0$;
\item $\Sigma$ is simply connected. 
\end{enumerate}
Then the map 
$$
H_{j}(\Sigma) \longrightarrow H_{j}(W)
$$
induced by inclusion is surjective.
\end{proposition}

\begin{proof} We can assume without loss of generality that the contact forms induced on $\Sigma$ viewed as contact hypersurface in $W$ and respectively in $P\times S^2$ are the same. Indeed, if $\alpha$ denotes the contact form coming from the contact embedding in $W$, and $\beta$ denotes the contact form coming from the contact embedding in $P\times S^2$, we can modify $W$ by attaching a large piece of symplectization $([1,R]\times \Sigma,d(r\alpha))$, $R\gg 1$ inside which we can find a graph over $\Sigma$ on which the induced contact form is a large multiple of $\beta$. By removing what lies beyond the graph and rescaling the symplectic form, we reduce ourselves to the situation where $\alpha=\beta$. 

We are now in a position to perform the construction of $V$ described above: we take away the interior $Z$ of $\Sigma$ in $P\times S^2$, and we replace it with $W$. 
Lemmas~\ref{lem:simple} and~\ref{lem:degree} hold true, and there exists an $\omega_V$-compatible regular almost complex structure $J$ on $V$ which satisfies the assumptions of Lemma~\ref{lem:degree}. The outcome is the compact manifold $\cM=\cM_J$ together with the degree $\pm 1$ map $\ev:\cM\times S^2\to V$.

Given the sets $S\subset T$, let us denote by $i^S_T$ the inclusion map of $S$ into $T$. Let $C_{1}\in H_{j}(W)$ be fixed. In order to prove that $C_1$ lies in the image of $(i^\Sigma_W)_*$, it is enough to show that  $C:=(i^W_V)_*(C_1)\in H_j(V)$ lies in the image of $(i^Y_V)_*$. 
Indeed, the  Mayer-Vietoris exact sequence writes
$$ \longrightarrow H_{j}(\Sigma) \stackrel {(i^{\Sigma}_W)_{*}\oplus (i^{\Sigma}_Y)_{*}}\longrightarrow H_{j}(W)\oplus H_{j}(Y)\stackrel{(i^{W}_V)_{*}-(i^{Y}_V)_{*}} \longrightarrow H_{j}(V) \longrightarrow H_{j-1}(\Sigma) \longrightarrow $$
Thus, if there is $C_{2}\in H_{j}(Y)$ such that $C=(i^{Y}_V)_{*}(C_{2})$, then  
$(C_{1},C_{2})$ is in the kernel of $(i^{W}_V)_{*}-(i^{Y}_V)_{*}$, hence in the image of
$(i^{\Sigma}_W)_{*}\oplus (i^{\Sigma}_Y)_{*}$. In particular $C_1$ is in the image of $(i^\Sigma_W)_*$.

We now prove that $C$ lies in the image of $(i_Y)_*$. We start with the observation that the map 
$$
(\ev)_{*}: H_{j}({\mathcal M} \times S^2) \longrightarrow H_{j}(V)
$$
is surjective. Since we use field coefficients, this is equivalent to injectivity of $\ev^*$ in cohomology, which in turn follows from the non-degeneracy of the cup-product pairing and the fact that $\ev:\mathcal M \times S^2 \longrightarrow V$ has non-zero degree (with respect to any field of coefficients). We can therefore write $C=\ev_*(\Gamma_C)$ for some $\Gamma_C\in H_j(\cM\times S^2)$ or, equivalently, 
$$
C=\ev_*(A\otimes \{pt\}+B \otimes [S^2])
$$ 
for some $A\in H_{j}(  {\mathcal M})$ and  $B \in H_{j-2}(  {\mathcal M})$. 

We claim 
that $B$ must vanish. Arguing by contradiction, let $B'$ be Poincar\'e dual to $B$ in $H_{*}( {\mathcal M})$, so that $B\cdot B'=\{pt\}$. We obtain  
 $$\Gamma_C\cdot (B'\otimes\{pt\}) = 
 (B\otimes [S^2]) \cdot (B'\otimes \{pt\}) = (B\cdot B') \otimes \{pt\}= \{pt\} \otimes \{pt\}.$$
This implies that
$$\{pt\}=(ev)_{*}(\Gamma_{C}\cdot (B'\otimes \{pt\}))= C \cdot (ev)_{*}(B'\otimes \{pt\}) = C \cdot ev^\infty_*(B'),$$
where $ev^z(u)=u(z)$. Since $ev^{\infty}_*(B') \subset P\times \{\infty\}$, we get $C \cdot ev^\infty_*(B')=0$, a contradiction.

As a result, we obtain $C=ev^z_*(A)$ (for any $z\in S^2$), with $A\in H_j(\mathcal M)$.
Choosing $z=\infty$ we get that $C \in (i^{P\times \{\infty\}}_V)_*(H_j(P)) \subset (i^Y_V)_*(H_j(Y))$. This concludes the proof. 

\end{proof}

\begin{remark} If the image of the boundary map $H_{j+1}(Y,\Sigma)\to H_j(\Sigma)$ coincides with the image of the boundary map 
$H_{j+1}(W,\Sigma)\to H_j(\Sigma)$, then $\dim H_{j}(W) \leq \dim H_{j}(P)$. Indeed, it follows from the commutative diagram below that 
the map $H_j(W)\to H_j(V)$ is injective. Since its image is contained in $\Image (i^{P\times\{\infty\}}_V)_*$, the conclusion follows. 
$$
\xymatrix
{
H_{j+1}(V,W) \ar[r] & H_j(W) \ar[r] & H_j(V) \\
H_{j+1}(Y,\Sigma) \ar[u]_{excision}^\simeq \ar[r]^\partial & H_j(\Sigma) \ar[r] \ar[u] & H_j(Y) \\
& H_{j+1}(W,\Sigma) \ar[u]_\partial & }
$$
\end{remark}

\begin{proof} [Proof of Theorem~\ref{prop:steinsubcritfillings}]
We use a result of Cieliebak (see \cite{Cieliebak1}) stating that a subcritical Stein manifold is symplectomorphic to $N\times {\mathbb C}$ where $N$ is Stein, and a result of Lisca and  Mati\'c (\cite{L-M}, Section 3, Theorem 3.2), stating that any Stein domain embeds symplectically  in a smooth projective manifold $P$ with ample canonical bundle. Moreover $N$ is contained in the complement of a hyperplane section $H$, which is of course hyperplane-section-like for the underlying complex structure.
Up to shrinking $\Sigma$ via the Liouville flow, we can thus assume that we have a contact embedding $\Sigma\subset N\times D^2\subset P\times S^2$, where $P$ carries an integral symplectic form $\omega_{P}$, the symplectic form $\sigma$ on $S^2$ is normalized by $[\sigma][S^2]=1$, and the image of $\Sigma$ is contained in $(P\setminus H)\times S^2$, where $H$ is a hyperplane-section-like symplectic submanifold. 
We may now apply Proposition~\ref{thm:BigMac} and this concludes our proof. 
\end{proof}

\begin{remark}
Of course the condition that $H_{j}(\Sigma) \longrightarrow H_{j}(W)$ is onto is equivalent to the claim that $H^j(W) \longrightarrow H^j(\Sigma)$ is injective, or that
$H_{j}(W) \longrightarrow H_{j}(W, \Sigma)$ vanishes.
\end{remark}

The case when $\Sigma$ is a sphere leads to the following variant of the Eliashberg-Floer-McDuff theorem (\cite{McDuff}): the assumptions that we impose are weaker, but so is the conclusion.
 \begin{Cor}
Let $(\Sigma, \xi)$ be a simply connected contact manifold admitting an embedding in a subcritical Stein manifold, and assume that $\Sigma$ is a homology sphere (resp. rational homology sphere). Any symplectically aspherical filling of $\Sigma$ is then a homology ball (resp. rational homology ball). 
\end{Cor}

\begin{proof}
Indeed apply Theorem \ref{prop:steinsubcritfillings} to the case where $H_{j}(\Sigma)=0$. We conclude that $H_{j}(W)=0$.
\end{proof}

(Counter-)examples are given by Brieskorn spheres (see Corollary \ref{Brieskorn}). 
Note that if $(\Sigma, \xi)$ is the standard contact sphere, it has an obvious embedding in $ {\mathbb R}^{2n}$. In this situation, using an argument by Eliashberg, it is proved in \cite{McDuff} that $W$ is simply connected. Thus we get, using  Smale's $h$-cobordism theorem (\cite{Smale}) that $W$ is diffeomorphic to the ball. This is the original Eliashberg-Floer-McDuff theorem.

 \begin{remark} 
 Here is a more precise statement.
 Remember that, in the proof of Theorem \ref{prop:steinsubcritfillings}, we showed that the image of $(i^{W}_V)_{*}$ is contained in the image of $(i^{P}_V)_{*}:=(i^{P\times \{\infty\}}_V)_{*}$ in $H_{*}(V)$. Now the following commutative diagram
 $$
\xymatrix
@C=40pt
{
H_j(W) \ar[r]^-{(i^W_V)_*} & H_j(V) & H_{j}(P) \ar[l]_-{(i^P_V)_*} \\
 H_j(\Sigma) \ar[r]^-{(i^\Sigma_Y)_*} \ar@{->>}[u] & H_j(Y)\ar[u] & H_{j}(P) \ar[l]_-{(i^P_Y)_*} \ar@{=}[u]  
}
$$
shows that $$\dim ((i^{P}_V)_{*}H_{j}(P)/(i^{W}_V)_{*}H_{j}(W)) \leq \dim  ((i^{P}_Y)_{*}(H_{j}(P)/(i^{\Sigma}_Y)_{*}H_{j}(\Sigma)).$$
We used that, given the linear map $\varphi=(i^Y_V)_*:H_j(Y)\to H_j(V)$, and linear subspaces $B\subseteq A\subseteq H_j(Y)$, we have 
$\dim\, \varphi(A)/\dim\, \varphi(B)\le \dim\, A/\dim\, B$.
 \end{remark}

\section{The case of $  ({\mathbb R} ^{2n}, \sigma_{0})$}

In this section we denote $b_{p}(X)$ the Betti numbers of a manifold $X$ with 
coefficients in a given field. Thus $b_p(X)$ is the rank of the $p$-th homology/cohomology group.

\begin{theorem} \label{thm:R2n}
Assume $(\Sigma, \xi)$ admits a contact embedding in $({\mathbb R} ^{2n}, \sigma_{0})$, with interior component $Z$.  
Let $(W,\omega)$ be a symplectic filling 
of $(\Sigma,\xi)$ such that $(\R^{2n}\setminus Z)\sqcup_\Sigma W$ is symplectically aspherical. This holds in particular if $W$ is symplectically aspherical and one of the following conditions is satisfied:
\begin{enumerate} 
\item \label{R:a} $H_{2}(W,\Sigma)=0$.
\item \label{R:b} The maps $\pi_{1}(\Sigma) \longrightarrow \pi_{1}(W)$ and $\pi_{1}(\Sigma) \longrightarrow \pi_{1}({\mathbb R}^{2n}\setminus Z)$ 
are injective. 
\item \label{R:c} $(\Sigma, \xi)$ is of restricted contact type in $(W, \omega)$ and in $({\mathbb R}^{2n}, \sigma_{0})$.
\end{enumerate} 

Then
\begin{enumerate}
\item[(1)] \label{item:1} any two symplectically aspherical fillings of $(\Sigma, \xi)$ which satisfy either of the conditions \ref{R:a}--\ref{R:c} 
have the same Betti numbers.
\item[(2)] \label{item:2} given a symplectically aspherical filling $W$ which satisfies one of the conditions \ref{R:a}--\ref{R:c}, 
the inclusion of $\Sigma$ in $W$ induces an injection in cohomology $$H^p(W) \longrightarrow H^p(\Sigma).$$
Moreover, we have $$b_{p}(\Sigma)= b_{p}(W)+b_{2n-p-1}(W).$$
\end{enumerate}
\end{theorem}
\begin{remark} Condition \ref{R:a} holds if $W$ is Stein and $n\ge 3$. Condition \ref{R:b} holds if $\Sigma$ is simply connected. 
The embedding $\Sigma\hookrightarrow \R^{2n}$ is always separating since $H_{2n-1}(\R^{2n};\Z)=0$.
\end{remark} 
\begin{remark}
The first statement in (2) follows from the previous section if either $\Sigma$ is simply connected or condition \ref{R:a} is satisfied. The reason why we can allow more general assumptions in the case of $\R^{2n}$ is that the geometry at infinity is perfectly controlled, unlike for an arbitrary subcritical Stein manifold. 
\end{remark}
\begin{remark} \label{rmk:symplectic}
There is a natural way to endow the smooth manifold $U:=({\mathbb R}^{2n} \setminus Z)\sqcup_\Sigma W$ with a symplectic form, which coincides with $\sigma_0$ on $\R^{2n}\setminus Z$. The assumption that $U$ is symplectically aspherical in the statement of Theorem \ref{thm:R2n} is understood with respect to this symplectic form. The construction is the following.  Let $\alpha_0$ be the contact form induced on $\Sigma$ as the concave boundary of $\R^{2n}\setminus Z$, and let $\alpha_W$ be the contact form induced on $\Sigma$ as the convex boundary of $W$. If $\alpha_0=\alpha_W$, then $\sigma_0$ and $\omega$ can be glued into a symplectic form on $U$. If $\alpha_0= f\alpha_W$ for some function $f:\Sigma\to(0,\infty)$, we reduce to the case of equality as follows. Let $m=\max_\Sigma (1/f)$ and choose $R> m\max_\Sigma(f)$. We attach to $W$ along its boundary the finite piece of symplectization $([1,R]\times \Sigma,d(r\alpha_W))$, remove what lies beyond the graph of $mf$, and denote the resulting domain by $W'$. Then $W'$ is diffeomorphic to $W$, it carries a natural symplectic form $\omega'$, while its contact boundary is naturally identified with $\Sigma$ and carries the contact form $mf\alpha_W=m\alpha_0$. Up to replacing 
$(W,\omega)$ with $(W',\frac 1 m \omega')$, we can therefore assume that $\alpha_W=\alpha_0$, and the two symplectic forms on $W$ and ${\mathbb R}^{2n} \setminus Z$ can be glued into a symplectic form on $U$. 
\end{remark}
\begin{proof}
We first show that any of the conditions \ref{R:a}--\ref{R:c} guarantees that the symplectic form $\omega'$ on  
$U:=({\mathbb R}^{2n} \setminus Z)\sqcup_\Sigma W$ described in Remark~\ref{rmk:symplectic} vanishes on spheres.  Since $\omega'=\sigma_0$ outside a compact set, we obtain that $U$ is diffeomorphic to ${\mathbb R}^{2n}$ by the Eliashberg-Floer-McDuff theorem. (Although for the convenience of the formulation we use the diffeomorphism statement in the Eliashberg-Floer-McDuff theorem, we only need the fact that $U$ has the homology of a point.)

Let $C$ be a $2$-sphere in $U$, assume without loss of generality that it intersects $\Sigma$ transversally, and denote by $C_1$ and $C_2$ the pieces contained in $W$ and $\R^{2n}\setminus \Sigma$ respectively. 
\begin{itemize}
\item Let us assume \ref{R:a}. Then we find a cycle $\Gamma$ in $\Sigma$ such that $\partial C_1=\partial \Gamma$. Since $C_1\cup \Gamma$ is a cycle in $H_2(W)$ and the map $H_2(\Sigma)\to H_2(W)$ is onto, we obtain that $C_1\cup\Gamma$ is homologous to a cycle $C_3$ contained in $\Sigma$. Since $\omega$ is exact near $\Sigma$, the area of $C_3$ is zero and therefore the areas of $C_1$ and of $\Gamma$ are equal. Hence $\Gamma\cup C_2$ is a cycle in $\R^{2n}\setminus Z$ with the same area as $C$. But the area of $\Gamma\cup C_2$ is zero because $\sigma_0$ is exact, and so is the area of $C$. 
\item Let us assume \ref{R:b}. At least one of the components of $C_1$ or $C_2$ is a disc, with boundary on $\Sigma$. By assumption, we can cap it by a disc in $\Sigma$ to get a sphere in $W$ or in $\R^{2n}\setminus Z$ which, by symplectic asphericity of $W$ and $\R^{2n}\setminus Z$, has zero area. We can thus inductively remove each component of $C_1$, $C_2$ and finally prove that $C$ has zero area. 
\item Let us assume \ref{R:c}. In this case the symplectic form on $U$ is exact. 
\end{itemize}

We thus have that $U$ is diffeomorphic to $\R^{2n}$. Since $\Sigma$ is contained in some large ball, denoted $B$, we equivalently have that $(B\setminus Z) \sqcup_\Sigma W$ is diffeomorphic to $B$.  The cohomology Mayer-Vietoris exact sequence writes in this case
$$
\longrightarrow H^p(B) \longrightarrow H^p(W)\oplus H^p(B\setminus Z)
\longrightarrow H^{p}(\Sigma) \longrightarrow H^{p+1}(B)
\longrightarrow
$$

Since $H^p(B)=0$ for $p > 0$, we see that the map $$H^p(W)\oplus H^p(B\setminus Z)
\longrightarrow H^{p}(\Sigma)$$ is an isomorphism for $p\geq 1$. Since it is induced by the inclusion maps,
the first claim in~(2) follows.

For  $p>0$   we have $$b_p(\Sigma)=b_p(W)+b_p(B\setminus Z).$$
Moreover, according to Alexander duality (see \cite{Greenberg-Harper}, theorem 27.5, p.233) we have $$b_p(B\setminus Z)=b_{2n-p-1}(Z)$$
for $0<p<2n-1$, which implies that, in this range, we have
$$b_p(\Sigma)=b_p(W)+b_{2n-p-1}(Z).$$
Of course, this also holds when we replace $W$ by $Z$, so that
$$b_p(\Sigma)-b_p(Z)=b_{2n-p-1}(Z)$$
and finally 
$$b_p(\Sigma)=b_p(W)+(b_p(\Sigma)-b_p(Z)).$$
This implies $b_p(W)=b_p(Z)$ for $0<p<2n-1$.

For $p=2n-1$,  if $B ( \varepsilon )$ is a small ball inside $Z$, the inclusions $$B\setminus B( \varepsilon ) \supset B\setminus  Z \supset S^{2n-1}$$
imply that  $b_{2n-1}(B-Z)\geq 1$, and the exact sequence
$$0 \longrightarrow H^{2n-1}(W)\oplus H^{2n-1}(B\setminus Z)
\longrightarrow H^{2n-1}(\Sigma) \longrightarrow 0$$
implies that $b_{2n-1}(B\setminus Z)=1$ and $b_{2n-1}(W)=b_{2n-1}(Z)=0$.
Finally, it is clear that the equality still holds for $p=0$, since $b_{0}(\Sigma)=b_{0}(W)=1$.
\end{proof}

\begin{Cor}\label{Cor3.5}
Assume $(\Sigma, \xi)$ has a Stein filling $(W, \omega)$ and has a contact
embedding in $ ({\mathbb R} ^{2n}, \sigma_{0})$, $n\ge 3$. Then
$$
\left\{ \begin{array}{l} b_{p}(\Sigma)=b_{p}(W)\ \text{for}\ 
0\leq p \leq n-2, \\ b_{n-1}(\Sigma)=b_{n}(\Sigma)=b_{n}(W)+b_{n-1}(W).\end{array}\right.
$$
Thus the homology of $W$ is completely determined by the homology of
$\Sigma$ except, maybe, in degree $n-1$ and $n$. It is completely determined by the homology of
$\Sigma$ if $b_{n}(\Sigma)=0$ or $W$ is subcritical Stein.  
\end{Cor}
\begin{proof}
The assumptions of Theorem~\ref{thm:R2n} are satisfied, since $H_2(W,\Sigma)=H^{2n-2}(W)$ vanishes if $n\ge 3$. It follows that $b_{p}(W)$  is determined by $b_{p}(\Sigma)$, except maybe in dimensions $n-1$ and $n$. If $b_{n}(\Sigma)=0$ we obtain
  $b_{n}(W)=b_{n-1}(W)=0$, and if $W$ is subcritical we have $b_{n}(W)=0$ and therefore $b_{n-1}(\Sigma)=b_{n-1}(W)$.
\end{proof}

\begin{remark}
Mei-Lin Yau proved  (see \cite{MLY}) that, if $W$ is subcritical Stein and the first Chern class of the complex vector bundle defined by $\xi$
vanishes, then cylindrical contact homology in the trivial homotopy class $HC_{*}^0(\Sigma, \alpha )$ is well-defined for a suitably chosen contact form $\alpha$, and we have
$$
HC_{*}^0(\Sigma, \alpha ) \simeq H_{*}(W,\Sigma)\otimes H_{*}( {\mathbb C}P^\infty)[2].
$$
By definition, the chain complex underlying $HC_*^0(\Sigma,\alpha)$ is generated by contractible closed Reeb orbits for the contact form $\alpha$. The degree of a generator $\gamma$ is defined to be $CZ(\gamma)+n-3$, where $CZ(\gamma)$ denotes the Conley-Zehnder index of the linearized Reeb flow along $\gamma$ in the transverse direction, computed with respect to a trivialization of $\xi$ along $\gamma$ induced by a trivialization over a spanning disc. The symbol $[2]$ denotes a shift in degree by $2$. 

It is therefore a general fact that the homology of a subcritical filling is determined by the contact structure
$(\Sigma, \xi)$. It is however not clear whether in general it is already determined by the knowledge
of the topology of $\Sigma$  (i.e. independently from $\xi$ or the topology of a filling).

\end{remark}

As a first consequence of Corollary~\ref{Cor3.5} and Mei-Lin Yau's result we have:

\begin{Cor}\label{3.7}
Assume $(\Sigma, \xi)$  satisfies $c_{1}(\xi)=0$, has a subcritical Stein filling $(W, \omega)$ and has a contact
embedding in $( {\mathbb R} ^{2n}, \sigma_{0})$. Then the rank of $HC_{*}^0(\Sigma, \alpha)$ is determined by $H_{*}(\Sigma)$.
Indeed, we have

$$
\rank(HC_{k}^0(\Sigma, \alpha))=
\sum_{\scriptsize \begin{array}{c} 2n-2-k \leq p\leq n-1 \\ p\equiv k\;  mod\; 2 \end{array} } 
b_{p}(\Sigma)
$$

\end{Cor}
\begin{proof} 
Note that assumption \ref{a}  from Theorem  \ref{thm:R2n} is automatically satisfied: we are in the Stein case. 
The result is a straightforward application of Corollary~\ref{Cor3.5}, Mei-Lin Yau's theorem and the duality $H^{2n-k}(W)\simeq H_{k}(W,\partial W)$. 

Thus $$HC_{k}^0(\Sigma,\alpha) = \bigoplus_{m\geq 0} H_{k-2m+2}(W,\Sigma) = \bigoplus_{m\geq 0} H^{2n-2-k+2m}(W) $$
and $b_{2n-2-k+2m}(W)=b_{2n-2-k+2m}(\Sigma)$ for $0\leq 2n-2-k+2m \leq n-1$. Setting $p=2n-2-k+2m$ yields the above formula. 
\end{proof} 
 
To state the next application of our theorem, let us recall the following definitions. A Hermitian line bundle $\mathcal L\stackrel \pi\to N$ over a symplectic manifold $(N^{2n-2}, \beta)$ is called {\it negative} if $c_1(\mathcal L)=-[\beta]$. Equivalently, there exists a Hermitian connection $\nabla$ whose curvature satisfies 
$\frac 1 {2i\pi}F^\nabla=-\beta$. Such a connection determines the {\it transgression $1$-form} $\theta^\nabla\in\Omega^1(\mathcal L\setminus 0_{\mathcal L},\R)$ which, by definition, vanishes on the horizontal distribution and is equal to $\frac 1 {2\pi}$ times the angular form in the fibers. Denoting $r(u):=|u|$, the total space $\mathcal L$ carries the symplectic form $\omega:=\pi^*\beta+d(r^2\theta^\nabla)$, which is exact on $\mathcal L\setminus 0_{\mathcal L}$ with $\omega=d((1+r^2)\theta^\nabla)$. The unit disc bundle $W=\{u\in\mathcal L\, : \, |u|\le 1\}$ is a symplectic manifold with contact type boundary. For details we refer to \cite{Oancea}, Section 3.3.

\begin{proposition} \label{3.8}
Let $(\Sigma,\xi)$ be the contact boundary of the unit disc bundle $(W,\omega)$ associated to a negative line bundle over $(N^{2n-2},\beta)$. Assume that $(N,\beta)$ is symplectically aspherical. Then  $(\Sigma, \xi )$ has no contact embedding in $({\mathbb R}^{2n}, \sigma_{0})$ with interior $Z$, such that $ ({\mathbb R} ^{2n}\setminus Z) \sqcup_\Sigma W$ is symplectically aspherical. The same holds for $n\geq 3$ and for any contact manifold obtained by contact surgery (as in \cite{Eliashberg, Weinstein}) of index  $k$ for any $k \in [3,n]$.
\end{proposition} 

\begin{proof}
The Gysin exact sequence writes 
$$ \longrightarrow H^{p-2}(N)\stackrel{ \beta \cup} \longrightarrow H^{p}(N) \longrightarrow H^p(\Sigma) \longrightarrow H^{p-1}(N) \longrightarrow $$
and, in degree $2$, we get
$$
H^{2}(\Sigma)=H^2(N)/ \langle [\beta]\rangle \oplus \ker \left ([\beta]\cup :H^{1}(N) \longrightarrow H^3(N)\right ).
$$
Hence $$b_{2}(\Sigma) < b_{2}(N)+b_{1}(N)=b_{2}(N)+b_{2n-2-1}(N)=b_2(W)+b_{2n-2-1}(W),$$ and this contradicts
Theorem \ref{thm:R2n}.

Let us now see what happens when we make a contact surgery. We shall denote our hypersurface by $\Sigma^-$, $W^-$ will be its filling,  and $\Sigma^+$ will be the result of the surgery on $\Sigma^-$ along a $(k-1)$-dimensional isotropic sphere.  Let us denote by $A_{k}\simeq D^k\times D^{2n-k}$ the attached handle, and denote $\partial ^-A_{k}= S^{k-1} \times D^{2n-k}$, $\partial ^+A_{k}= D^{k} \times S^{2n-k-1}$, 
so that 
the new filling of $\Sigma^+$  is $ W^+=W^- \cup_{\partial^- A_{k}} A_{k}$. 
We first need to prove that $W^+$ is symplectically aspherical. But the homotopy exact sequence of the pair $(W^+,W^-)$ is given by 
$$ \longrightarrow  \pi_{3}(W^+, W^-) \longrightarrow \pi_{2}(W^-) \longrightarrow \pi_{2}(W^+) \longrightarrow \pi_{2}(W^+, W^-) \longrightarrow $$
and $\pi_{2}(W^+, W^-)\simeq \pi_{2}(A_{k}, \partial^{-}A_{k}) \simeq \pi_{2}(D^{k}, \partial D^{k})= 0$ for $k\geq 3$. Thus the inclusion of $W^-$ in $W^+$ induces a surjective map on $\pi_{2}$, hence if $[\omega]\pi_{2}(W^-)=0$, we also have $[\omega] \pi_{2}(W^+)=0$.

Let us now  first consider the case $k\geq 4$. We claim that we have $b_{2}(\Sigma^+)= b_{2}(\Sigma^-)$ and $b_{2}(W^+)= b_{2}(W^-)$. Indeed the homology exact sequence for the pair $(W^+, W^-)$ writes
$$ \longrightarrow  H_{3}(W^+, W^-) \longrightarrow H_{2}(W^-) \longrightarrow H_{2}(W^+) \longrightarrow H_{2}(W^+, W^-) \longrightarrow $$
but $H_{j}(W^+, W^-)\simeq H_{j}(A_{k}, \partial^{-}A_{k}) \simeq H_{j}(D^{k}, \partial D^{k})= 0$ for $j=2,3$ and $k\geq 4$, so $b_{2}(W^+)=b_{2}(W^-)$. 

Similarly  the Mayer-Vietoris exact sequence for $\Sigma^\pm= \Sigma^-\setminus (\partial ^-A_{k}) \cup \partial^\pm A_{k}$ reads
  \begin{equation}\label{3.1} \begin{array}{ll} H_{2}(S^{k-1} \times S^{2n-k-1}) \longrightarrow H_{2}(\Sigma^- \setminus \partial^{-}A_{k}) \oplus H_{2}(\partial^\pm A_{k}) \ \\ \hspace{5cm} \longrightarrow H_{2}(\Sigma^\pm) \longrightarrow H_{1}(S^{k-1} \times S^{2n-k-1}) = 0  \end{array} \end{equation}    
When  $k\geq 4$, the groups  $H_{2}(S^{k-1} \times S^{2n-k-1})$ and $H_{2}(\partial^\pm A_{k})$ vanish, so that we have isomorphisms
$$ H_{2}(\Sigma^- \setminus \partial^{-}A_{k}) \simeq H_{2}(\Sigma^\pm)$$ and therefore $b_{2}(\Sigma^+)= b_{2}(\Sigma^-)$.

\medskip 

Let us now deal with the case $k=3, n\geq 4$. 
 In case ``$-$'', the first map in the exact sequence \eqref{3.1} is injective (since its projection on the second summand is induced by the inclusion $S^{2}\times S^{2n-4} \longrightarrow S^{2}\times D^{2n-3}$ ). Since $2n-4 >2$, we obtain $b_{2}(\Sigma^-)=b_{2}(\Sigma^-\setminus \partial^-A_{3})$. In the  ``$+$'' case, we have  
$H_{2}(\partial^+A_{3})=0$ and $b_{2}(\Sigma^+)\leq b_{2}(\Sigma^-\setminus \partial^-A_{3})=b_{2}(\Sigma^-)$.

We write the homology exact sequences of the pairs $(W^+,W^-)$ and $(\Sigma^+,\Sigma^+\cap W^-)$

{\footnotesize
$$
\xymatrix 
@C=10pt
{H_{3}(W^+,W^-) \ar[r]^{\partial_{W}} & H_{2}(W^-) \ar[r] &H_{2}(W^+) \ar[r] & H_{2}(W^+,W^-) =0 \\ 
H_{3}(\Sigma^+, \Sigma^+\cap W^-) \ar[r]^{\partial_{\Sigma}} \ar[u]& H_{2}(\Sigma^+\cap W^-) \ar[r] \ar[u]&H_{2}(\Sigma^+) \ar[r] \ar[u]& H_{2}(\Sigma^+, \Sigma^+\cap W^-)=0  \ar[u]& \\ & H_{2}(\Sigma^-\setminus \partial^-A_{3})\simeq H_{2}(\Sigma^-) \ar@{=}[u] & &
}$$
}

The left hand side vertical map is an isomorphism since 
\begin{gather*} H_{3}(\Sigma^+, \Sigma^+\cap W^-)\simeq H_{3}(D^3\times S^{2n-4}, S^2\times  S^{2n-4}) \hspace{4cm} \\  \hspace{5cm}  \stackrel{\simeq}\longrightarrow H_{3}(D^3\times D^{2n-3},S^2\times D^{2n-3}) \simeq H_{3}(W^+,W^-)  \end{gather*} 
Therefore either the map $\partial_{W}$ is injective, and thus so is $\partial_{\Sigma}$ and consequently $b_{2}(W^+)=b_{2}(W^-)-1$ and $b_{2}(\Sigma^+)= b_{2}(\Sigma^{-})-1$,  or it is zero, and then $b_{2}(W^+)=b_{2}(W^-)$ and $b_{2}(\Sigma^+)\leq  b_{2}(\Sigma^{-})$. 

For $n=k=3$, 
we leave it to the reader to check that  $$H_{2}(\Sigma^-)\simeq H_{2}(\Sigma^- \setminus \partial^-A_{3})/ \Image (H_{2}(\{pt\}\times S^2))$$
$$H_{2}(\Sigma^+)\simeq H_{2}(\Sigma^- \setminus \partial^-A_{3})/ \Image (H_{2}(S^2\times S^2))$$ so that $b_{2}(\Sigma^+)$ equals either $b_{2}(\Sigma^-)$ or $b_{2}(\Sigma^-)-1$. Again using the same argument as above, whenever   $b_{2}(W^+)= b_{2}(W^-)-1$ we have $b_{2}(\Sigma^+)=  b_{2}(\Sigma^-)-1$.  This concludes our proof. \end{proof}
\begin{remark} 
According to \cite{Laudenbach}, if $(\Sigma, \xi)$ has a contact embedding in $ {\mathbb R}^{2n}$, the same holds for any manifold obtained by contact surgery over an isotropic sphere of dimension $\leq n-1$. In contrast, we display here an obstruction to embedding $\Sigma$ in $ {\mathbb R}^{2n}$ that survives such a surgery.  
\end{remark}

\begin{examples}
The symplectic asphericity condition in Proposition~\ref{3.8} is necessary:  the manifold $({\mathbb C}P^{n-1}, \sigma_{0})$ is not symplectically aspherical, and $(S^{2n-1}, \alpha_{0})$ has a contact embedding into $( {\mathbb R}^{2n}, \sigma_{0})$.
 \end{examples}
 
\begin{remark} The previous proof does not generalize to higher rank bundles. Let us call a Hermitian vector bundle $\mathcal E\stackrel \pi\to N$ over a symplectic manifold $(N^{2n-2}, \beta)$  {\it negative} if it admits a Hermitian connexion $\nabla$ whose curvature $\frac 1 i F^\nabla\in\Omega^2(N,\mathrm{End}\, \mathcal E)$ is negative definite. This means that, for any $\beta$-compatible almost complex structure $J$ on the base $N$ and any non-zero vector $v\in TN$, we have $\frac 1 i F^\nabla(v,Jv)<0$. 

Let $\mathbb{P}(\mathcal E)$ denote the projectivized bundle and $\mathcal L_{\mathcal E}\stackrel \pi \to \mathbb{P}(\mathcal E)$ be the tautological line bundle. Then $\mathcal L_{\mathcal E}$ is a negative Hermitian line bundle, and the total space carries the symplectic form $\omega_{\mathcal L}=\pi^*\omega_{FS}+\Omega_{\mathcal L}$, where $\omega_{FS}$ is the curvature form on $\mathbb{P}(\mathcal E)$ and $\Omega_{\mathcal L}=d(r^2\theta^\nabla)$, with $r(u)=|u|$ and $\theta^\nabla$ the transgression $1$-form (see the preamble to Proposition~\ref{3.8}). 
Denoting $W_{\mathcal L}=\{u\in\mathcal L_{\mathcal E}\, : \, |u|\le 1\}$, we see that $\Sigma=\partial W_{\mathcal L}$ is a contact manifold. However, the filling $W_{\mathcal L}$ is \emph{not} symplectically aspherical since it contains $\mathbb{P}(\mathcal E)$ as a symplectic submanifold. 

The manifold $\Sigma$ can also be realized inside $\mathcal E$ as $\{u\in\mathcal E\, : \, |u|= 1\}$, via a natural diffeomorphism $\mathcal L_{\mathcal E}\setminus 0_{\mathcal L_{\mathcal E}}\simeq \mathcal E\setminus 0_{\mathcal E}$. This diffeomorphism transforms $W_{\mathcal L}$ into $W=\{u\in\mathcal E\, : \, |u|\le 1\}$. The pull-back of $\Omega_{\mathcal L}$ via this diffeomorphism, denoted $\Omega$, is symplectic on $\mathcal E\setminus 0_{\mathcal E}$ and extends over $0_{\mathcal E}$, as equal to the area form in the fibers and vanishing along the zero-section. We can thus equip $\mathcal E$ with the symplectic form $\omega=\pi^*\beta+\Omega$. If $\beta$ is symplectically aspherical, then so is $\omega$. However, $\Sigma=\partial W$ is \emph{not} of contact type since the restriction of $\pi^*\beta$ to $\Sigma$ is not exact for $r\ge 2$, as shown by the Gysin exact sequence. 

The outcome of this discussion is that, even if $(N,\beta)$ is symplectically aspherical, $\Sigma$ does not appear naturally as contact type boundary of a symplectically aspherical manifold. We feel that a result analogous to Proposition~\ref{3.8} should hold for higher rank negative vector bundles, but our methods do not apply in this case. 

For the details of the above constructions we refer to \cite{Oancea}, Section 3.4. 
\end{remark}

 \begin{proposition}\label{3.12}
 Let $L$ be a compact manifold admitting a Lagrangian embedding into $ {\mathbb R}^{2n}$ and $n\geq 3$.  Then any symplectically aspherical filling $W$ of
 $ST^*L$  such that $H_2(W,ST^*L)=0$ has the same homology as $DT^*L$ (and hence the homology of $L$).
 \end{proposition}
 \begin{proof}
 Indeed, the hypothesis implies that $ST^*L$ has a contact (non exact !) embedding into $ {\mathbb R} ^{2n}$, so that
 we can apply Theorem \ref{thm:R2n}. The condition $H_{2}(DT^*L,ST^*L)=0$ is clearly satisfied using Thom's isomorphism. 
 \end{proof}

 Let now $ST^*L$ be the unit cotangent bundle of $L$. The spectral sequence of this sphere bundle yields the following dichotomy:

\begin{itemize}\item either the Euler class vanishes, and then $$b_{p}(ST^*L)= b_{p}(L)+b_{p-(n-1)}(L)$$ 

\item or the Euler class is non zero and then

 $$\left\{ \begin{array}{ll} b_{p}(ST^*L)= b_{p}(L)+b_{p-(n-1)}(L) \;\text{for}\;
p\neq n-1,n \\ b_{n}(ST^*L)=b_{n-1}(ST^*L)= b_{n-1}(L)=b_{1}(L)\end{array}\right.$$
\end{itemize}

The formula $$b_{p}(\Sigma)=b_{p}(W)+b_{2n-p-1}(W)$$ becomes

\begin{enumerate}
\item in the first case

$$b_{p}(L)+b_{p-(n-1)}(L)=b_{p}(L)+b_{2n-p-1}(L)$$
hence $$b_{p-(n-1)}(L)=b_{2n-p-1}(L)$$ that is the Poincar\'e duality formula
\item
in the second case

$$b_{1}(L)=b_{n-1}(L)=b_{n}(L)+b_{2n-n-1}(L)=b_{n}(L)+b_{n-1}(L)$$
\end{enumerate}
This implies $b_{n}(L)=0$, which is impossible (at least for orientable $L$).

\begin{proposition}
Let $L$ be an orientable manifold with non zero Euler class. Then $ST^*L$ has no contact embedding in $ {\mathbb R}^{2n}$, $n\ge 3$. The same holds for any contact manifold obtained from such a $ST^*L$ by surgery of index $3\leq k\leq n-3$.
\end{proposition}
\begin{proof}
The case of $ST^*L$ has been already proved above. The surgery does not modify the conditions $H_{2}(W,\Sigma)=0$ nor does it change $b_{n}(\Sigma)$ or $b_{n}(W), b_{n-1}(W)$. This concludes our proof. 
\end{proof} 

\begin{remark}
The condition $e(L)=0$ is exactly the condition needed to be able to find a Lagrangian immersion of $L$ regularly homotopic to an embedding. We however suspect that there are no embeddings of $ST^*L$ as a a smooth hypersurface in $ {\mathbb R}^{2n}$. 
\end{remark}

 \section{The Stein subcritical case}

 In this section we assume that $(\Sigma_{1}, \xi_{1})$ has a separating contact embedding in a 
 subcritical Stein domain $(W_{2}, \omega_{2})$ with boundary $(\Sigma_{2}, \xi_{2})$, and we denote 
 by $V_1$ the bounded component of $W_2\setminus \Sigma_1$. We denote by $(W_{1}, \omega_{1})$ an arbitrary symplectically aspherical filling of $(\Sigma_1,\xi_1)$ such that one of the following assumptions holds (cf. Theorem~\ref{prop:steinsubcritfillings}) 
 \begin{itemize}
\item $H_{2}(W_1,\Sigma_1)=0$.
\item $\Sigma_1$ is simply connected.
 \end{itemize}
  
 \begin{proposition}
 Under the above assumptions, we have that
 $$b_{j}(W_{1})\leq b_{j}(\Sigma_{1})+ \min(0, b_{j}(\Sigma_{2})-b_{j}(W_{2}\setminus V_{1 }))$$ 
 \end{proposition}

\begin{proof}
Note that given an exact sequence $A \overset {f}\longrightarrow B \overset {g}\longrightarrow C$ we have
$\dim (B) = \dim (\ker (g))+ \dim (\Image (g)) = \dim (\Image f) + \dim (\Image (g)) \leq \dim (A) + \dim (C)$.
Using the Mayer-Vietoris exact sequence of $(W_{2}\setminus V_{1}) \sqcup W_{1}$ and the inequality 
$\dim H_{j}(\Sigma_{2}) \geq \dim H_{j}((W_{2}\setminus V_{1}) \sqcup W_{1})$, we get that
$$b_{j}(W_{2}\setminus V_{1}) + b_{j}(W_{1}) \leq b_{j}(\Sigma_{2}) + b_{j}(\Sigma_{1}).$$
Thus $$b_{j}(W_{1})\leq (b_{j}(\Sigma_{2})-b_{j}(W_{2}\setminus V_{1 }))+b_{j}(\Sigma_{1}).$$
Acccording to Theorem \ref{prop:steinsubcritfillings} we have  $b_{j}(W_{1})\leq b_{j}(\Sigma_{1})$,  and our claim follows. 
\end{proof}

 Note that $b_{j}(W_{2}\setminus V_{1})=b_{2n-j}(W_{2},V_{1} \cup \Sigma_2)$ by Poincar\'e duality and excision.
Note also that the above result is stronger than Theorem \ref{prop:steinsubcritfillings} only when $ b_{j}(\Sigma_{2})-b_{j}(W_{2}\setminus V_{1 })<0$. 
This happens for example if $\Sigma_{2}$ is a homology sphere.

 \medskip 
 
 The first part of the following result has been obtained in a weaker form and by different methods in \cite{Cieliebak-Frauenfelder-Oancea} (see also proposition \ref{5.11}). 
\begin{proposition} \label{4.4}
Let $L$ be an orientable closed manifold of dimension $\geq 3$, with non-zero Euler class. 
Then $ST^*L$ has no contact embedding in a subcritical Stein manifold. As before, this also holds for any manifold 
obtained from $ST^*L$ by contact surgery of index $k \in [3,n-1]$. 
\end{proposition} 
 \begin{proof} 
Since $n \geq 3$, the group $H_{2}(DT^*L, ST^*L)$ is zero, so assumption \ref{a} of Theorem \ref{prop:steinsubcritfillings} is satisfied. The Gysin exact sequence of $ST^*L$ shows that the map $H_{n}(ST^*L) \longrightarrow H_{n}(L)$ vanishes.   This contradicts Theorem \ref{prop:steinsubcritfillings}.
The case of manifolds obtained by surgery is dealt with as in Proposition \ref{3.8}.  
\end{proof}

 \section{Obstructions from Symplectic homology}\label{shof}

  In this section we assume that $(\Sigma, \xi)$ is a contact manifold whose first Chern class $c_{1}(\xi )$ vanishes. 
  All the symplectic fillings $(W,\omega)$ of $(\Sigma, \xi)$ that we consider are assumed to satisfy $c_{1}(TW)=0$.

 \begin{definition}
Let $(W, \omega)$ be a connected symplectically aspherical manifold with contact type boundary. We say that $(W, \omega)$ is an {\it $SAWC$ manifold} if  $SH_{*}(W, \omega)=0$.
\end{definition}

\begin{remark}
The vanishing of $SH_*(W)$ is equivalent to the fact that $W$ satisfies the {\it Strong Algebraic Weinstein Conjecture} as defined in \cite{Viterbo1}, stating that the canonical map $H_{2n}(W,\partial W)\to SH_n(W)$ is not injective. This follows from the following three observations: (i) non-injectivity of the map $H_{2n}(W,\partial W)\to SH_n(W)$ is equivalent to its vanishing, since $H_{2n}(W,\partial W)$ is $1$-dimensional; (ii) symplectic homology is a ring with unit \cite{McLean}, and the unit is the image of the fundamental class of $W$ under the map $H_{2n}(W,\partial W)\to SH_n(W)$ \cite{Seidel}; (iii) vanishing of the unit for $SH_*(W)$ is equivalent to the vanishing of $SH_*(W)$.  

It is proved in \cite{Bourgeois-Oancea-3}, Corollary~1.4 that an $SAWC$ manifold also satisfies  
the {\it Equivariant Algebraic Weinstein Conjecture} from \cite{Viterbo1}. This can also be seen using the spectral sequence connecting the usual version of symplectic homology to the equivariant version \cite{Viterbo1}.
\end{remark}

 If we have an exact embedding $(V_{1},\omega_{1})$ into $(W_{1}, \omega_{1})$,  there is an induced transfer  map (see \cite{Viterbo1})  $$SH_{*}(W_{1}) \longrightarrow SH_{*}(V_{1})$$ which, according to  Mark McLean (see \cite{McLean}), is a unital ring homomorphism. This implies the following result:
\begin{proposition} [\cite{McLean}]\label{prop:weakly exact}
Let $(V,\omega)$ be an exact symplectic submanifold of $(W, \omega)$. If $(W, \omega)$ is  $SAWC$ then $(V,\omega)$
is also $SAWC$.
\end{proposition}

It is easy to find $SAWC$ manifolds which are not Stein. For example, we have:

\begin{proposition}[\cite{Oancea}]
Let $P$ be any exact symplectic manifold with contact type boundary. Then, for any exact $SAWC$ manifold $W$, 
we have that $P\times W$ is
$SAWC$. Also, the total space of a negative symplectic fibration in the sense of \cite{Oancea} with fiber $W$ is $SAWC$.
\end{proposition}

\begin{proposition} [\cite{Cieliebak2}]
Let $W'$ be obtained from $W$ by attaching handles of index $\leq n-1$. Then $SH_{*}(W)\simeq SH_{*}(W')$. In particular if $W$ is $SAWC$, the same holds for $W'$.
\end{proposition} 

The following statement is contained 
in \cite{Cieliebak-Frauenfelder-Oancea}, Corollary~1.15 and Remark~1.19.

  \begin{theorem} \label{thm:subcrit}
 Let $(\Sigma, \xi)$ be a contact manifold for which there exists a 
 contact form $\alpha$ whose closed characteristics are nondegenerate and have Conley-Zehnder 
 index strictly bigger than $3-n$. Let $i:(\Sigma, \xi)\hookrightarrow (W,\omega)$ be a separating  
 exact embedding in an $SAWC$ manifold $(W, \omega)$. Assume $i_*:\pi_1(\Sigma)\to \pi_1(W)$ is injective.
  Then the Betti numbers of the interior $V$ of $\Sigma$ in any coefficient field are 
 determined by the contact structure $\xi$ (and do not depend on the choice of the $SAWC$ manifold $W$).
 \end{theorem}
 \begin{proof}
By Proposition~\ref{prop:weakly exact} we have  $SH_{*}(V)=0$. 
The relative exact sequence in symplectic homology (see \cite{Viterbo1}) then implies that

 $$SH^+_{*}(V)\simeq H_{*+n-1}(V, \partial V).$$

We can assume without loss of generality that the contact form induced by the contact embedding on $\Sigma$ is equal to $\alpha$. Let $\widehat V=(V,\omega)\cup([1,\infty[\times \Sigma,d(r\alpha))$ be the symplectization of $V$, obtained by gluing a semi-infinite cone along the boundary. 
If the Reeb vector field associated to the contact form $\alpha$ has no closed characteristic with Conley-Zehnder index $\le 3-n$, and 
 if $i_*:\pi_1(\Sigma)\to \pi_1(V)$ is injective, then  
 there is no rigid holomorphic plane in $\widehat V$ bounding a closed characteristic. 
In this case, it is a consequence of the stretch-of-the-neck argument in \cite{Bourgeois-Oancea-2} that $SH_*^+(V)$ depends only on the contact boundary $\partial V=\Sigma$ (see also \cite{Cieliebak-Frauenfelder-Oancea}, Corollary~1.15).  As a consequence, the Betti numbers $b_{j}(V, \partial V)$ only depend on $\xi$.

 Now if we have another exact embedding of $\Sigma$ in $W'$ and $W'$ is also $SAWC$, the interior $V'$  of $\Sigma$ in 
 $W'$ must have the same cohomology as $V$.
 \end{proof}

 \begin{proposition} \label{5.6}
 Let $(\Sigma,\xi)$ be the boundary of a subcritical Stein manifold $(W,\omega)$. Let $(M,\omega)$ be an $SAWC$
 manifold such that $(\Sigma, \xi)$ has an exact separating embedding into $(M,\omega)$, with interior $Z$. 
 Then $H_{*}(Z)\simeq H_{*}(W)$. 
 \end{proposition} 
 \begin{proof} 
 First of all, by Proposition \ref{prop:weakly exact}, we have $SH_{*}(Z)=0$. On one hand the exact sequence (\cite{Viterbo1})
 $$\longrightarrow SH_{*}(Z) \longrightarrow SH_{*}^+(Z)  \longrightarrow H_{*+n-1}(Z,\Sigma) \longrightarrow SH_{*-1}(Z) \longrightarrow $$
 shows that $H_{*}(Z,\Sigma)\simeq SH_{*+1-n}^+(Z)$. On the other hand, since $\Sigma$ bounds a subcritical Stein manifold, there exists a contact form 
 $\alpha$ such that the Reeb orbits are all non-degenerate and of index $>3-n$ (cf. \cite{MLY}), so the proof of Theorem~\ref{thm:subcrit} implies that $SH^+_{*}(Z)\simeq SH^+_{*}(W)$. This last space is in turn isomorphic to 
 $H_{*+n-1}(Z,\Sigma)$ by the same argument, and finally $H_{*}(Z,\Sigma) \simeq H_{*}(W,\Sigma)$, hence $H_{*}(Z)\simeq H_{*}(W)$. 
 \end{proof}  
 
 \begin{remark} 
 The condition that $W$ is subcritical is not really necessary. We only need $W$ to be $SAWC$ provided  there is a contact form defining  $\xi$ for which all closed Reeb orbits are nondegenerate and have index  $>  3-n$.
 \end{remark} 
 \begin{remark} Proposition~\ref{5.6} can be compared to the following result: 
 \begin{Cor}(\cite{MLY}) \label{cor:subcrit}
 Let $W$ be a subcritical Stein manifold with boundary $\partial W$ such that $c_1(TW)|_{\pi_2(W)}=0$. 
 Then any subcritical Stein manifold with the same  boundary $\partial W$ and whose first Chern class vanishes on the second homotopy group 
 has the same homology as $W$.
 \end{Cor}

 \begin{proof}

This follows from  the main computation of \cite{MLY}
$$HC_{*}(\p W, \alpha ) \simeq H_{*}(W,\p W)\otimes H_{*}( {\mathbb C}P^\infty),$$
which implies directly that the homology of a subcritical Stein filling is determined by the contact structure of the boundary.
 \end{proof}
 \end{remark} 
 \begin{remark}
 Note that, when $\Sigma=S^{2n-1}$, we may apply our proposition to $W=D^{2n}$. Thus we prove that any symplectically aspherical filling $Z$ with vanishing first Chern class satisfies $H_{*}(Z)=0$ in nonzero  degree. Thus, if $Z$ is simply connected and $n\geq 3$, it is diffeomorphic to a ball. This is a weak version of the Eliashberg-Floer-McDuff theorem mentioned in the previous section, but note that the above proof does not make use of it and also that it extends to many other contact manifolds.

 \end{remark}

Let us now use the above tools to find obstructions to contact embeddings. 
We first have: 

\begin{proposition} [\cite{Cieliebak-Frauenfelder-Oancea}] \label{5.11}
  If $(\Sigma,\xi) = (ST^*L,\xi_{\mathrm{std}})$ with  $L$ a closed simply connected manifold, then
  $(\Sigma, \xi)$  has no separating exact embedding in an $SAWC$ manifold $(M, \omega)$. Here $\xi_{\mathrm{std}}$ denotes the standard contact structure on $ST^*L$.
 \end{proposition} 
 \begin{proof} Since the characteristic flow on $ST^*L$ is the geodesic flow, it has all closed trajectories of index $\geq 0 >3-n$ if $n>3$ (in the cases $n=2,3$ we have that  $L$ is a sphere and we can find a metric for which all closed geodesics have index $>3-n$). 
Assuming the existence of such an embedding, with interior $Z$, the proof of Theorem~\ref{thm:subcrit} shows that $SH_*^+(Z)$ depends only on the boundary $(\Sigma,\xi)$. We obtain on the one hand $SH_*^+(Z)\simeq H_{*+n-1}(Z,\partial Z)$, and on the other hand $SH_*^+(Z)\simeq SH_*^+(DT^*L)$. But $SH_*^+(DT^*L)\simeq H_*(\Lambda L,L)$, where $\Lambda L$ denotes the free loop space of $L$. Hence $SH_*^+(DT^*L)$ is infinite dimensional, a contradiction. 
 \end{proof} 
 
 \begin{remark} Let $(M, \omega)$ be obtained by attaching subcritical handles to $DT^*L$. Provided one can prove that the Reeb orbits on $(\partial M, \xi)$ still have index $>3-n$, our argument extends to show that $(\partial M, \xi)$ has no contact embedding in an $SAWC$ manifold.
 \end{remark}

The case of circle bundles can also be dealt with using  contact and Floer homology, as follows.

\begin{proposition}\label{5.13}
 Let $(\Sigma, \xi)$ be the unit circle bundle associated to a negative line bundle $\mathcal L$ over a symplectically aspherical manifold $(N^{2n-2r}, \beta)$ such that  $c_{1}(TN)=0$. Then, for $n\geq 2$,   $\Sigma$ does not bound a subcritical Stein manifold with vanishing first Chern class. 
 The same holds for any contact manifold obtained by subcritical surgery on $(\Sigma, \xi)$ of index $\neq 2, 3$. 
\end{proposition}

\begin{proof}
Indeed, let $M$ denote the manifold bounded by $\Sigma$. If $W$ is the unit disc bundle associated to $\Sigma$, we have $\partial W= \Sigma$ and,
using that $SH_{*}(W)=0$ (\cite{Oancea}) and the exact sequence

$$  \longrightarrow SH_{*}(W) \longrightarrow SH^+_{*}(\Sigma) \longrightarrow H_{*+n-1}(W,\Sigma) \longrightarrow $$
we obtain   $$SH^+_{*}(\Sigma) \simeq H_{*+n-1}(W,\Sigma) \simeq H_{*+n-3}(N).$$
The same exact sequence with $M$ yields $$SH^+_{*}(\Sigma) \simeq H_{*+n-1}(M,\Sigma) \simeq H^{n-*+1}(M)$$
But this last space vanishes for $*\leq 1$ while $H_{*+n-3}(N)$ is non-zero for $*=3-n$. When $n\geq 2$ we get a contradiction. 
Now since $k\neq 2, 3$, $H_{2}(W,\Sigma)$ does not change, so remains equal to $H_{0}(N)= {\mathbb Q}$. But we must have $H_{2}(W,\Sigma)=SH_{3-n}^+(\Sigma)=H^{2n-2}(M)=0$.  A contradiction. \end{proof}

\begin{remark}
This partially answers a question of Biran in \cite{Biran} who asked the same question in the Stein case (not subcritical). A different answer was given by
Popescu-Pampu in \cite{Popescu-Pampu}
\end{remark}

\section{Brieskorn manifolds, McLean's examples }
We consider an isolated  singularity of holomorphic germ. For example, assume we are given $V$ a complex submanifold in $ {\mathbb C}^{n+1}$ 
with an isolated singularity at the origin. We then consider the
submanifold $\Sigma_{ \varepsilon }= S_{ \varepsilon } \cap V$, where $S_{ \varepsilon } =\{z \in {\mathbb C}^{n+1} \mid \vert z \vert ^2= \varepsilon \}$ and $\varepsilon>0$ is small enough.
The maximal complex subspace of the tangent space defines a hyperplane distribution which happens to be a contact structure, 
and whose isotopy class is independent of $ \varepsilon $. In case the singularity is smoothable, $\Sigma_{ \varepsilon }$ bounds a Stein manifold $W$.

\begin{example}
If $V$ is the hypersurface $f^{-1}(0)$ where $f$ is polynomial, then the singularity is always smoothable.  
According to \cite{Milnor}, the manifold $W$ is homotopy equivalent to a wedge of $n$-spheres. The number $\mu\ge 0$ of spheres is called {\em the Milnor number of the singularity}. We obtain that $W$ is $(n-1)$-connected and $H_n(W)=  {\mathbb Z}^{\mu}$. 

The boundary $\Sigma:=\partial W$ is $(n-2)$-connected. It is called the {\it link of the singularity}. For $n\ge 2$ the long exact sequence of the pair $(W,\Sigma)$ reduces to 
\begin{equation} \label{eq:seqMilnor}
0\to H_n(\Sigma)\to H_n(W)\stackrel {S} \to \mathrm{Hom}(H_n(W),\Z)\to H_{n-1}(\Sigma)\to 0. 
\end{equation}
Here we used the identification $H_n(W,\Sigma)\simeq H^n(W)\simeq \mathrm{Hom}(H_n(W),\Z)$. It is proved in \cite{Milnor} that the map $S$ is given by the intersection form, namely  
$$
S(x)(y):=x\cdot y, \qquad x,y\in H_n(W).
$$

One also defines the {\it Seifert form of the singularity} 
$$
A:H_n(W)\otimes H_n(W)\to \Z
$$ 
by $A(x,y):=\mathrm{lk}_{S_\varepsilon}(x^+,y)$, where $W$ is now viewed inside $S_\varepsilon$ via the Milnor open book given by $f/|f|$, $y^+$ denotes a small push-off of $y$ in the positive direction given by the open book decomposition, and $\mathrm{lk}_{S_\varepsilon}$ denotes the linking number inside $S_\varepsilon$. We then have 
$S=A+(-1)^n {A^t}$ (see~\cite{Durfee} and the references therein). 
\end{example}

As an immediate consequence of Proposition \ref{prop:steinsubcritfillings}, we obtain the following result.

 \begin{proposition}\label{prop:Brieskorn}
(a)  Let $n\geq 3$ and  $(\Sigma, \xi)$ be the link of an isolated hypersurface singularity. If the intersection form on the middle-dimensional homology of the Milnor fiber is nonzero, then $(\Sigma,\xi)$ does not embed in a subcritical Stein manifold. 

(b) Brieskorn manifolds of dimension $2n-1$, $n\ge 3$ with Milnor number at least $2$ do not admit contact embeddings in subcritical Stein manifolds. 
 \end{proposition}

 \begin{proof}
 (a) The long exact sequence~\eqref{eq:seqMilnor} shows that surjectivity of the map $H_n(\Sigma)\to H_n(W)$ is equivalent to the vanishing of the intersection form. Since $n\ge 3$ we have $H_2(W,\Sigma)=0$, and the conclusion follows from Theorem~\ref{prop:steinsubcritfillings}.
 
 (b) The Brieskorn manifold $\Sigma(a_0,a_1,\dots,a_n)$, $a_0,\dots,a_n\ge 2$ is, by definition, the link of the singularity $z_0^{a_0}+\dots +z_n^{a_n}=0$. The Milnor number of  $\Sigma(a_0,a_1,\dots,a_n)$ is $\mu=(a_0-1)\dots(a_n-1)$. Following \cite{Sakamoto}, its Seifert form is the tensor-product of blocks of dimension $a_i-1$, $i=0,\dots,n$, and the blocks have the form \cite{Durfee} 
 $$
 \left(\begin{array}{ccccc}
 1 & 1 & 0 & \dots & 0 \\
 0 & 1 & 1 & \ddots & 0 \\
 0 & \ddots & 1 & \ddots  & 0 \\
 \vdots &          &    & \ddots & 1 \\
  0 & & \dots & 0 & 1
 \end{array}\right).
 $$
 Thus $A$ is neither symmetric, nor anti-symmetric, and we infer that $S\neq 0$. The conclusion then follows from (a).
 \end{proof}
 
 \begin{remark} The condition $\mu\ge 2$ is violated if and only if all the exponents $a_i$ are equal to $2$. In this case $\Sigma=ST^*S^n$. For $n$ even the matrix of $A$ is symmetric, hence $S=A+A^t\neq 0$, so there is no contact embedding of $\Sigma$ in a subcritical Stein manifold. But an argument of Lisca in~\cite{Cieliebak-Frauenfelder} shows that there is not even a smooth embedding. If $n$ is odd we cannot conclude.  
 \end{remark}
 
 \begin{Cor} \label{Brieskorn} Let $n\ge 3$ and $(\Sigma,\xi)$ be a Brieskorn manifold which is diffeomorphic to the sphere $S^{2n-1}$. The standard contact structure $\xi$ inherited from the Milnor fiber is exotic, i.e. $\xi$ is not diffeomorphic to the standard contact structure on $S^{2n-1}$.
 \end{Cor}
 
 \begin{proof}
 This follows immediately from Proposition~\ref{prop:Brieskorn}, since $\Sigma$ does not admit a contact embedding in $\R^{2n}$. There is no need to consider the case $a_i=2$ for all $i$ since $ST^*S^n$ is never diffeomorphic to $S^{2n-1}$. 
 \end{proof}

 \begin{remark}
Ustilovsky has actually exhibited in \cite{Usti} infinitely many pairwise non-isomorphic contact structures on spheres of dimension $4m+1$. In \cite{vanKoert}, the reader will  find an algorithm to compute the linearized contact homology of most Brieskorn manifolds  in dimension greater than $5$. 
 \end{remark}
 
Let us now consider the manifolds of Mark McLean in \cite{McLean}. These are Stein symplectic manifolds $(M^{2n}_{k}, \omega_{k})$ diffeomorphic to $ {\mathbb R} ^{2n}$ ($n\geq 4$), such that $(\partial M_{k}^{2n}, \xi_{k}^{2n})$ is a contact manifold diffeomorphic to $S^{2n-1}$. However, McLean shows that $SH_{n}(M_{k}^n)$ contains $N^k$ idempotent elements for some $N\ge 2$, therefore the manifolds $M_{k}^{2n}$ are pairwise non symplectomorphic. 

We now prove

\begin{proposition} \label{6.3}
The contact manifolds $(\partial M_{k}^{2n},\xi_{k}^{2n})$ are never contactomorphic to the standard sphere. 
\end{proposition} 
\begin{proof} Let us denote for simplicity $W= M_{k}^{2n}$ and $(\Sigma, \xi) =(\partial M_{k}^{2n},\xi_{k}^{2n})$. The exact sequence in symplectic homology reads
$$ \longrightarrow H_{2n}(W,\Sigma) \longrightarrow SH_{n}(W) \longrightarrow SH_{n}^+(W) \longrightarrow 0$$
Assume $(\Sigma, \xi)$ is the standard sphere. Then $SH_{n}^+(W)$ only depends on $(\Sigma, \xi)$ so is the same as $SH_{n}^+(D^{2n})=0$. 
As a result we should have $\rank (SH_{n}(W)) \leq 1$. But for $k\geq 2$, there are at least $3$ idempotents, hence the  rank is at least $2$ and 
we get a contradiction. 
\end{proof} 

If we knew that there is a contact form on $(\partial M_{k}^{2n},\xi_{k}^{2n})$ with no closed characterisitic of index less than $3-n$, then we would get, by the above argument,  that $(\partial M_{k}^{2n},\xi_{k}^{2n})$ has no embedding in an $SAWC$ manifold. 
 \section{Summary}

A conceptual framework for the study of symplectic fillings is provided by the following definition of \cite{E-H}. 

\begin{definition}[\cite{E-H}] \label{def:dominated} Let $(\Sigma_{1}, \alpha_{1})$ and $ (\Sigma_{2}, \alpha_{2})$ be two closed contact manifolds. We say that $(\Sigma_{1}, \alpha_{1})$ 
is {\bf dominated by} $ (\Sigma_{2}, \alpha_{2})$  if there exists a symplectically aspherical manifold $(W, \omega)$ such that $(W, \omega)$ has  $(\Sigma_{1}, \alpha_{1})$ as a concave boundary, $ (\Sigma_{2}, \alpha_{2})$ as a convex boundary and no other boundary component. We shall write $$(\Sigma_{1}, \alpha_{1}) \prec  (\Sigma_{2}, \alpha_{2}).$$ We shall say that $(\Sigma_{1}, \alpha_{1})$ is {\bf equivalent to}
$(\Sigma_{2}, \alpha_{2})$ if we both have $(\Sigma_{1}, \alpha_{1}) \prec  (\Sigma_{2}, \alpha_{2})$ and $(\Sigma_{2}, \alpha_{2}) \prec  (\Sigma_{1}, \alpha_{1})$, and this is denoted by
$$(\Sigma_{1}, \alpha_{1})\simeq (\Sigma_{2}, \alpha_{2}).$$
\end{definition}

\begin{remark}
In the terminology of Symplectic Field Theory, we see that $(\Sigma_{1}, \alpha_{1})$ 
is dominated by $ (\Sigma_{2}, \alpha_{2})$ if and only if there exists a symplectically aspherical cobordism between $(\Sigma_{1}, \alpha_{1})$ and $ (\Sigma_{2}, \alpha_{2})$. 
\end{remark}

Clearly, we have
$$(\Sigma_{1}, \alpha_{1})\simeq (\Sigma_{1}, \alpha_{1}).$$
We would like to know if there are nonequivalent pairs of contact manifolds. Clearly, a contact manifold admits a filling if and only 
if it dominates the standard sphere. Which manifolds are dominated by the standard sphere? Our results give examples of 
fillable manifolds which are not dominated by the standard sphere or, more generally, by the boundary of a subcritical Stein manifold.
On the other hand, in dimension $4$, any overtwisted contact manifold is dominated by any other contact manifold (see \cite{E-H}). 
In particular, all overtwisted contact structures are equivalent! 

The point of view of Definition~\ref{def:dominated} is also related to the work of \cite{Chantraine} on the non-symmetry of Legendrian concordances. 
 \newpage
We here try to summarize our results, but warn the reader that in the table below,  the assumptions of the theorems are usually incomplete and the statements often not precise. One should refer to the relevant section of the paper for full details. 
\label{table}\begin{center} \footnotesize
    \begin{xtabular}{| p{2.2cm} | p{3.8cm} | p{5 cm}  | p{4cm} |} \hline
  & \bf Weakly subcritical case&   \bf Stein subcritical case & \bf Case of $ {\mathbb R}^{2n}$  \\ \hline
    {\it Hypothesis A} & $(\Sigma, \alpha)$ has a separating contact embedding   in  an $SAWC$ manifold $(M,\omega)$ with bounded component $Z$.    & $(\Sigma, \alpha)$ has a contact embedding   in a subcritical Stein
     $(M ,\omega)$ with bounded component $Z$.   &   $(\Sigma, \alpha)$ has a contact embedding   in $ {\mathbb R}^{2n}$ with bounded component $Z$  \\  \hline
    {\it Assume} & $(\Sigma, \alpha)=\partial (W, \omega_{1})$ & $(\Sigma, \alpha)=\partial (W, \omega_{1})$   &   $(\Sigma, \alpha)=\partial (W, \omega_{1})$ \\ \hline
    {\it Conclusion 1} & &    The map $H_{j}(\Sigma) \longrightarrow H_{j}(W)$ is onto (Thm. \ref{prop:steinsubcritfillings}) & The homology of $W$ is (almost)  determined by the homology of $\Sigma$ (Thm. \ref{thm:R2n}) \\ \hline      {\it Hypothesis B} & $W$ is subcritical Stein &   $W$ is $SAWC$ & $W$ is subcritical Stein  \\ \hline
    {\it Conclusion} &  $(\Sigma, \xi)$  determines the homology of $Z$ (Prop.  \ref{5.6})  and the rank of $HC_{*}(\Sigma)$ is determined by $H_{*}(W)$ (\cite{MLY})  &  If the Conley-Zehnder indices of closed characteristics are $> 3-n$,  $(\Sigma, \xi)$  determines the homology of $Z$ (Prop.  \ref{5.6})  and the rank of $HC_{*}(\Sigma)$ is determined by $H_{*}(W)$ (\cite{MLY}) &  The rank of $HC_{*}(\Sigma)$ is determined by $H_{*}(\Sigma)$  (Prop. \ref{3.7}) \\ \hline
    {\it Examples: uniqueness of fillings}  &  & Any filling of a simply connected homology sphere embeddable in a subcritical Stein is a homology ball.  &If $L$ has a Lagrange embedding in $ {\mathbb R}^{2n}$, the fillings of $ST^*L$ have the homology of $L$.    (Prop. \ref{3.12}) \\   \hline    
     {\it Examples: obstructions to contact embeddings} &  & - Circle bundles of negative line bundles and some of   their surgeries having no contact embedding in a subcritical Stein (Prop. \ref{5.13}) &  Circle bundles of negative line bundles and some of  their surgeries having no contact embedding in $ {\mathbb R}^{2n}$ (Prop. \ref{3.8})  \\  
         &  Obstructions to contact embedding  $ST^*L$  in an $SAWC$ manifold.  (Prop. \ref{5.11}) & - Obstructions to contact embedding  $ST^*L$ and the manifolds obtained from it by surgery in a subcritical Stein.    (Prop. \ref{4.4}) \newline - Brieskorn manifolds do  not embed in subcritical Stein (Prop.~\ref{prop:Brieskorn}) \newline - Brieskorn spheres are exotic contact spheres (Cor. \ref{Brieskorn})  \newline - Contact spheres obtained by \cite{McLean} as boundaries of exotic symplectic $ {\mathbb R}^{2n}$  are exotic  (Cor. \ref{6.3}) &\\   \hline 
     \end{xtabular}
\end{center}

\end{document}